\documentclass{amsart}

\newtheorem{theorem}{Theorem}[section]

\theoremstyle{definition}

\theoremstyle{remark}

\numberwithin{equation}{section}

\begin{document}
\title[Identities for generalized twisted Bernoulli polynomials
]{$\begin{array}{c}
         \text{Identities of symmetry for generalized twisted Bernoulli polynomials}\\
         \text{twisted by ramified roots of unity}\
       \end{array}$}

\author{dae san kim} \thanks{}
\address{Department of Mathematics, Sogang University, Seoul 121-742, Korea}
\email{dskim@sogong.ac.kr}
\address{}
\email{}
\thanks{}

\subjclass[]{}

\date{}

\dedicatory{}


\begin{abstract}
We derive eight identities of symmetry in three  variables related to generalized twisted Bernoulli polynomials and generalized twisted power sums, both of which are twisted by ramified roots of unity. All of these are new, since there have been results only about identities of symmetry in two variables. The derivations of identities are based on the $p$-adic integral expression of the generating function for the generalized twisted Bernoulli polynomials and the quotient of
$p$-adic integrals that can be expressed as the exponential generating function for the generalized twisted power sums.

\end{abstract}

\maketitle

\section{Introduction and preliminaries}
 Let $p$ be a fixed prime. Throughout this paper,
$\mathbb{Z}_{p}$, $\mathbb{Q}_{p}$, $\mathbb{C}_{p}$ will respectively denote the ring of
$p$-adic integers, the field of $p$-adic rational numbers and the completion of the algebraic closure of
$\mathbb{Q}_{p}$. Assume that $|~|_{p}$ is the normalized absolute value of $\mathbb{C}_{p}$ such that
$|p|_{p}=\frac{1}{p}$.  The group $\Gamma$ of all roots of unity of
$\mathbb{C}_{p}$ is the direct product of its subgroups
$\Gamma_{u}$(the subgroup of unramified roots of unity) and
$\Gamma_{r}$(the subgroup of ramified roots of unity). Namely,

\begin{align*}
\Gamma=\Gamma_{u}\cdot\Gamma_{r},~
\Gamma_{u}\cap\Gamma_{r}=\{1\},
\end{align*}
where
\begin{align*}
\Gamma_{u}=\{\xi\in\mathbb{C}_{p}|\xi^{r}=1~for~
some~r\in\mathbb{Z}_{>0}~with~(r,p)=1\},\\
\Gamma_{r}=\{\xi\in\mathbb{C}_{p}|\xi^{p^{s}}=1
~for~some~s\in\mathbb{Z}_{>0}\}.
\end{align*}
Let $d$ be a fixed positive integer. Then we let
\begin{align*}
X=X_{d}=\lim_{\overleftarrow{N}}\mathbb{Z}/dp^{N}\mathbb{Z},
\end{align*}
and let $\pi:X\rightarrow\mathbb{Z}_{p}$ be the map given by the inverse limit of the natural maps
\begin{align*}
\mathbb{Z}/dp^{N}\mathbb{Z}\rightarrow\mathbb{Z}/p^{N}\mathbb{Z}.
\end{align*}

If g is a function on $\mathbb{Z}_{p}$,we will use the same notation to denote the function
$g\circ\pi$.  Let $\chi:(\mathbb{Z}/d\mathbb{Z})^{*}\rightarrow
\overline{\mathbb{Q}}^{*}$ be a (primitive) Dirichlet character of conductor $d$.  Then it will be pulled back to
$X$ via the natural map $X\rightarrow\mathbb{Z}/d\mathbb{Z}$.
Here we fix, once and for all, an imbedding
$\overline{\mathbb{Q}}\rightarrow\mathbb{C}_{p}$, so that
$\chi$ is regarded as a map of $X$ to $\mathbb{C}_{p}$(cf.\cite{R8})

For a uniformly differentiable function
$f:X\rightarrow\mathbb{C}_{p}$, the $p$-adic Volkenborn-type integral of
$f$ is defined(cf.\cite{R5}) by
\begin{align*}
\int_{X}f(z)d\mu(z)=\lim_{N\rightarrow
\infty}\frac{1}{dP^{N}}\sum_{j=0}^{dp^{N}-1}f(j).
\end{align*}
Then it is easy to see that

\begin{equation}\label{a1}
\int_{X}f(z+1)d\mu(z)=\int_{X}f(z)d\mu(z)
+f'(0).
\end{equation}
More generally, we deduce from (\ref{a1}) that, for any positive integer $n$,

\begin{equation}\label{a2}
\int_{X}f(z+n)d\mu(z)=\int_{X}f(z)d\mu(z)
+\sum_{a=0}^{n-1}f'(a).
\end{equation}
Throughout this paper, we let
$\xi\in\Gamma_{r}$ be any fixed root of unity, and let
\begin{equation}\label{a3}
E=\{t\in\mathbb{C}_{p}||t|_p< p^{\frac{-1}{p-1}}\}.
\end{equation}
Then, for each fixed $t\in E$, the function
$e^{zt}$ is analytic on $\mathbb{Z}_{p}$ and hence considered as a function on $X$ and, by applying (\ref{a2}) to $f$ with
$f(z)=\chi(z)\xi^{z}e^{zt}$, we get the $p$-adic integral expression of the generating function for the generalized twisted Bernoulli numbers $B_{n,\chi,\xi}$ attached to $\chi$ and $\xi$;

\begin{equation}\label{a4}
\int_{X}\chi(z)\xi^{z}e^{zt}d\mu(z)
=\frac{t}{\xi^{d}e^{dt}-1}\sum_{a=0}^{d-1}\chi(a)\xi^{a}e^{at}
=\sum^{\infty}_{n=0}B_{n,\chi,\xi}\frac{t^{n}}{n!}(t\in E).
\end{equation}

So we have the following $p$-adic integral
expression of the generating function for the
generalized twisted Bernoulli polynomials
$B_{n,\chi.\xi}(x)$ attached to $\chi$ and $\xi$;

\begin{equation}\label{a5}
\int_{X}\chi(z)\xi^{z}e^{(x+z)t}d\mu(z)
=\frac{te^{xt}}{\xi^{d}e^{dt}-1}\sum_{a=0}^{d-1}\chi(a)\xi^{a}e^{at}
=\sum^{\infty}_{n=0}B_{n,\chi,\xi}(x)\frac{t^{n}}{n!}(t\in E, x\in\mathbb{Z}_{p}).
\end{equation}
Also, from (\ref{a1}) we have:
\begin{equation}\label{a6}
\int_{X}\xi^{z}e^{zt}d\mu(z)=\frac{t}{\xi e^{t}-1}.
\end{equation}

 Let $S_{k}(n;\chi,\xi)$ denote the $k$th generalized twisted power sum of the first $n+1$ nonnegative integers attached to
$\chi$ and $\xi$, namely

\begin{equation}\label{a7}
S_{k}(n;\chi,\xi)=\sum_{a=0}^{n}\chi(a)\xi^{a}a^{k}
=\chi(0)\xi^{0}0^{k}+\chi(1)\xi^{1}1^{k}+\cdots +\chi(n)\xi^{n}n^{k}.
\end{equation}
From (\ref{a4}), (\ref{a6}), and (\ref{a7}), one easily derives the following identities: for $w\in\mathbb{Z}_{>0}$,

\begin{align}
\label{a8}
\frac{dw\int_{X}\chi(x)\xi^{x}e^{xt}d\mu(x)}{\int_{X}
\xi^{dwy}e^{dwyt}d\mu(y)}
&=\frac{\xi^{dw}e^{dwt}-1}{\xi^{d}e^{dt}-1}
\sum_{a=0}^{d-1}\chi(a)\xi^{a}e^{at}\\
\label{a9}
&=\sum_{a=0}^{dw-1}\chi(a)\xi^{a}e^{at}\\
\label{a10}
&=\sum_{k=0}^{\infty}S_{k}(dw-1;\chi,\xi)\frac{t^{k}}{k!}~
(t\in E).
\end{align}
In what follows, we will always assume that the $p$-adic integrals of the various (twisted) exponential functions on $X$ are defined for
$t\in E$(cf. (\ref{a3})), and therefore it will not be mentioned.

  \cite{R1}, \cite{R2}, \cite{R6}, \cite{R9} and \cite{R10} are some of the previous works on identities of symmetry in two variables involving Bernoulli polynomials and power sums. For the brief history, one is referred to those papers. For the first time, the idea of \cite{R6} was extended in \cite{R4} to  the case of three variables so as to yield many new identities with abundant symmetry. This added some new identities of symmetry even to the existing ones in two variables as well.

  In this paper, we will produce  8 basic identities of symmetry in three variables $w_1$,$w_2$,$w_3$ related to generalized twisted
 Bernoulli polynomials and generalized twisted
 power sums, both of which are twisted by ramified roots of unity(i.e., $p$-power roots of unity)(cf.
 (\ref{a38})-(\ref{a41}), (\ref{a44})-(\ref{a47}).
 All of these seem to be new,  since there have been results only about identities of symmetry in two variables in the literature
(\cite{R7}).  On the other hand, in \cite{R3} the measure introduced by Koblitz(cf.\cite{R8})was exploited in order to treat the unramified roots of unity case (i.e., the orders of the roots of unity are prime to $p$ and the conductors of Dirichlet characters).

The following is stated as Theorem \ref{A5} and  an example of the full six symmetries in the positive integers $w_1$,$w_2$,$w_3$.

\begin{equation*}
\begin{split}
w_{1}^{n-1}\sum_{k=0}^{n}\binom{n}{k}\sum_{a=0}^{dw_{1}-1}
\chi(a)\xi^{aw_{2}w_{3}}B_{k,\chi,\xi^{w_{1}w_{3}}}
&(w_{2}y_{1}+\frac{w_{2}}{w_{1}}a)\\
&\times S_{n-k}(dw_{3}-1;\chi,\xi^{w_{1}w_{2}})
w_{2}^{n-k}w_{3}^{k-1}
\end{split}
\end{equation*}
\begin{equation*}
\begin{split}
=w_{1}^{n-1}\sum_{k=0}^{n}\binom{n}{k}\sum_{a=0}^{dw_{1}-1}
\chi(a)\xi^{aw_{2}w_{3}}B_{k,\chi,\xi^{w_{1}w_{2}}}
&(w_{3}y_{1}+\frac{w_{3}}{w_{1}}a)\\
&\times S_{n-k}(dw_{2}-1;\chi,\xi^{w_{1}w_{3}})
w_{3}^{n-k}w_{2}^{k-1}
\end{split}
\end{equation*}
\begin{equation*}
\begin{split}
=w_{2}^{n-1}\sum_{k=0}^{n}\binom{n}{k}\sum_{a=0}^{dw_{2}-1}
\chi(a)\xi^{aw_{1}w_{3}}B_{k,\chi,\xi^{w_{2}w_{3}}}
&(w_{1}y_{1}+\frac{w_{1}}{w_{2}}a)\\
&\times S_{n-k}(dw_{3}-1;\chi,\xi^{w_{1}w_{2}})
w_{1}^{n-k}w_{3}^{k-1}
\end{split}
\end{equation*}
\begin{equation*}
\begin{split}
=w_{2}^{n-1}\sum_{k=0}^{n}\binom{n}{k}\sum_{a=0}^{dw_{2}-1}
\chi(a)\xi^{aw_{1}w_{3}}B_{k,\chi,\xi^{w_{1}w_{2}}}
&(w_{3}y_{1}+\frac{w_{3}}{w_{2}}a)\\
&\times S_{n-k}(dw_{1}-1;\chi,\xi^{w_{2}w_{3}})
w_{3}^{n-k}w_{1}^{k-1}
\end{split}
\end{equation*}
\begin{equation*}
\begin{split}
=w_{3}^{n-1}\sum_{k=0}^{n}\binom{n}{k}\sum_{a=0}^{dw_{3}-1}
\chi(a)\xi^{aw_{1}w_{2}}B_{k,\chi,\xi^{w_{2}w_{3}}}
&(w_{1}y_{1}+\frac{w_{1}}{w_{3}}a)\\
&\times S_{n-k}(dw_{2}-1;\chi,\xi^{w_{1}w_{3}})
w_{1}^{n-k}w_{2}^{k-1}
\end{split}
\end{equation*}
\begin{equation*}
\begin{split}
=w_{3}^{n-1}\sum_{k=0}^{n}\binom{n}{k}\sum_{a=0}^{dw_{3}-1}
\chi(a)\xi^{aw_{1}w_{2}}B_{k,\chi,\xi^{w_{1}w_{3}}}
&(w_{2}y_{1}+\frac{w_{2}}{w_{3}}a)\\
&\times S_{n-k}(dw_{1}-1;\chi,\xi^{w_{2}w_{3}})
w_{2}^{n-k}w_{1}^{k-1}.
\end{split}
\end{equation*}

The derivations of identities are based on the $p$-adic integral
expression of the generalized twisted Bernoulli polynomials in
 (\ref{a5}) and the quotient of integrals in (\ref{a8})-(\ref{a10}) that can be
expressed as the exponential generating function for the generalized twisted power sums.  These abundance of symmetries would not be unearthed if such $p$-adic integral representations had not been available. We indebted this idea to the paper \cite{R7}.

\section{Several types of quotients of $p$-adic integrals}

  Here we will introduce several types of quotients of $p$-adic integrals on $X$ or $X^{3}$ from which some
interesting identities follow owing to the built-in symmetries in
$w_{1},w_{2},w_{3}$. In the following, $w_{1},w_{2},w_{3}$ are all
positive integers and all of the explicit expressions of integrals
in (\ref{a12}), (\ref{a14}), (\ref{a16}), and (\ref{a18}) are obtained from the identity in (\ref{a4}) and (\ref{a6}). To ease notations, from now on we suppress
$\mu$ and denote, for example, $d\mu(x)$ and $d\mu(x_{1})d\mu(x_{2})d\mu(x_{3})$ respectively simply by
$dx$ and $dX$.
\\
\\
(a) Type $\Lambda_{23}^{i}$ (for $i=0,1,2,3$)\\
$I(\Lambda_{23}^{i})$
\begin{equation}
\label{a11}
\quad=\frac{\begin{array}{c}
         d^{i}\int_{X^{3}}\chi(x_{1})\chi(x_{2})\chi(x_{3})
\xi^{w_{2}w_{3}x_{1}+w_{1}w_{3}x_{2}+w_{1}w_{2}x_{3}} \\
         \qquad\qquad\qquad\qquad\qquad\times e^{(w_{2}w_{3}x_{1}+w_{1}w_{3}x_{2}+w_{1}w_{2}x_{3}
+w_{1}w_{2}w_{3}(\sum_{j=1}^{3-i}y_{j}))t}dX
       \end{array}
}
{(\int_{X}\xi^{dw_{1}w_{2}w_{3}x_{4}}
e^{dw_{1}w_{2}w_{3}x_{4}t}dx_{4})^{i}}
\end{equation}
\begin{align}
\begin{split}
\label{a12}
=&\frac{(w_{1}w_{2}w_{3})^{2-i}t^{3-i}
e^{w_{1}w_{2}w_{3}
(\sum_{j=1}^{3-i}y_{j})t}(\xi^{dw_{1}w_{2}w_{3}}
e^{dw_{1}w_{2}w_{3}t}-1)^{i}}
{(\xi^{dw_{2}w_{3}}e^{dw_{2}w_{3}t}-1)
(\xi^{dw_{1}w_{3}}e^{dw_{1}w_{3}t}-1)
(\xi^{dw_{1}w_{2}}e^{dw_{1}w_{2}t}-1)}\\
&\times(\sum_{a=0}^{d-1}\chi(a)\xi^{aw_{2}w_{3}}e^{aw_{2}w_{3}t})
(\sum_{a=0}^{d-1}\chi(a)\xi^{aw_{1}w_{3}}e^{aw_{1}w_{3}t})
(\sum_{a=0}^{d-1}\chi(a)\xi^{aw_{1}w_{2}}e^{aw_{1}w_{2}t}).
\end{split}
\end{align}
\\
\\
(b) Type $\Lambda_{13}^{i}$ (for $i=0,1,2,3$)\\
$I(\Lambda_{13}^{i})$
\begin{align}
\label{a13}
&=\frac{d^{i}\int_{X^{3}}\chi(x_{1})\chi(x_{2})
\chi(x_{3})\xi^{w_{1}x_{1}+w_{2}x_{2}+w_{3}x_{3}}
e^{(w_{1}x_{1}+w_{2}x_{2}+w_{3}x_{3}+w_{1}w_{2}w_{3}
(\sum_{j=1}^{3-i}y_{j}))t}dX}
{(\int_{X}\xi^{dw_{1}w_{2}w_{3}x_{4}}e^{dw_{1}w_{2}w_{3}x_{4}t}dx_{4})^{i}}\\
\label{a14}
&=\frac{(w_{1}w_{2}w_{3})^{1-i}t^{3-i}e^{w_{1}w_{2}w_{3}
(\sum_{j=1}^{3-i}y_{j})t}(\xi^{dw_{1}w_{2}w_{3}}
e^{dw_{1}w_{2}w_{3}t}-1)^{i}}{(\xi^{dw_{1}}e^{dw_{1}t}-1)
(\xi^{dw_{2}}e^{dw_{2}t}-1)(\xi^{dw_{3}e^{dw_{3}t}}-1)}
\end{align}
\begin{align*}
\qquad \qquad \qquad \qquad\qquad \times (\sum_{a=0}^{d-1}\chi(a)\xi^{aw_{1}}e^{aw_{1}t})
(\sum_{a=0}^{d-1}\chi(a)\xi^{aw_{2}}e^{aw_{2}t})
(\sum_{a=0}^{d-1}\chi(a)\xi^{aw_{3}}e^{aw_{3}t}).
\end{align*}
\\
\\
(c-0) Type $\Lambda_{12}^{0}$\\
$I(\Lambda_{12}^{0})$
\begin{align}
\label{a15}
&=\int_{X^{3}}\chi(x_{1})\chi(x_{2})\chi(x_{3})
\xi^{w_{1}x_{1}+w_{2}x_{2}+w_{3}x_{3}}
e^{(w_{1}x_{1}+w_{2}x_{2}+w_{3}x_{3}+w_{2}w_{3}y+w_{1}w_{3}y+w_{1}
w_{2}y)t}dX
\\
\label{a16}
&=\frac{w_{1}w_{2}w_{3}t^{3}e^{(w_{2}w_{3}+w_{1}w_{3}+w_{1}w_{2})yt}}
{(\xi^{dw_{1}}e^{dw_{1}t}-1)
(\xi^{dw_{2}}e^{dw_{2}t}-1)(\xi^{dw_{3}}e^{dw_{3}t}-1)}
\end{align}
\begin{align*}
\qquad\times (\sum_{a=0}^{d-1}\chi(a)\xi^{aw_{1}}e^{aw_{1}t})
(\sum_{a=0}^{d-1}\chi(a)\xi^{aw_{2}}e^{aw_{2}t})
(\sum_{a=0}^{d-1}\chi(a)\xi^{aw_{3}}e^{aw_{3}t}).
\end{align*}
\\
\\
(c-1) Type $\Lambda_{12}^{1}$
\begin{align}
\label{a17}
I(\Lambda_{12}^{1})&=\frac{d^{3}\int_{X^{3}}\chi(x_{1})\chi(x_{2})
\chi(x_{3})\xi^{w_{1}x_{1}+w_{2}x_{2}+w_{3}x_{3}}
e^{(w_{1}x_{1}+w_{2}x_{2}+w_{3}x_{3})t}dX}{\int_{X^{3}}
\xi^{d(w_{2}w_{3}z_{1}+w_{1}w_{3}z_{2}+w_{1}w_{2}z_{3})}
e^{d(w_{2}w_{3}z_{1}+w_{1}w_{3}z_{2}+w_{1}w_{2}z_{3})t}dZ}\\
\label{a18}
&=\frac{(\xi^{dw_{2}w_{3}}e^{dw_{2}w_{3}t}-1)
(\xi^{dw_{1}w_{3}}e^{dw_{1}w_{3}t}-1)
(\xi^{dw_{1}w_{2}}e^{dw_{1}w_{2}t}-1)}
{w_{1}w_{2}w_{3}(\xi^{dw_{1}}e^{dw_{1}t}-1)
(\xi^{dw_{2}}e^{dw_{2}t}-1)(\xi^{dw_{3}}e^{dw_{3}t}-1)}
\end{align}
\begin{align*}
\qquad \qquad \qquad\times (\sum_{a=0}^{d-1}\chi(a)\xi^{aw_{1}}e^{aw_{1}t})
(\sum_{a=0}^{d-1}\chi(a)\xi^{aw_{2}}e^{aw_{2}t})
(\sum_{a=0}^{d-1}\chi(a)\xi^{aw_{3}}e^{aw_{3}t}).
\end{align*}
All of the above $p$-adic integrals of various types are invariant
under all permutations of $w_{1},w_{2},w_{3}$, as one can see either
from $p$-adic integral representations in (\ref{a11}), (\ref{a13}),
(\ref{a15}), and (\ref{a17}) or from their explicit evaluations in
(\ref{a12}), (\ref{a14}), (\ref{a16}), and (\ref{a18}).

\section{Identities for generalized twisted Bernoulli polynomials}
All of the following results can be easily obtained from (\ref{a5})
and (\ref{a8})-(\ref{a10}). First, let's consider Type
$\Lambda_{23}^{i}$, for each $i=0,1,2,3.$
\\
\\
(a-0)
\begin{equation*}
\begin{split}
&I(\Lambda_{23}^{0})\\
&=\int_{X}\chi(x_{1})\xi^{w_{2}w_{3}x_{1}}e^{w_{2}w_{3}(x_{1}+w_{1}y_{1})
t}dx_{1}\int_{X}\chi(x_{2})\xi^{w_{1}w_{3}x_{2}}
e^{w_{1}w_{3}(x_{2}+w_{2}y_{2})t}dx_{2}
\qquad\qquad\qquad\qquad\qquad\qquad\\
\end{split}
\end{equation*}
\begin{equation*}
\qquad\qquad\qquad\qquad\qquad\qquad\qquad\qquad\qquad\times
\int_{X}\chi(x_{3})\xi^{w_{1}w_{2}x_{3}}
e^{w_{1}w_{2}(x_{3}+w_{3}y_{3})t}dx_{3}
\end{equation*}
\begin{equation*}
=(\sum_{k=0}^{\infty}\frac{B_{k,\chi,\xi^{w_{2}w_{3}}}(w_{1}y_{1})}
{k!}(w_{2}w_{3}t)^{k})(\sum_{l=0}^{\infty}\frac{B_{l,\chi,\xi^{w_{1}w_{3}}}
(w_{2}y_{2})}{l!}(w_{1}w_{3}t)^{l})\qquad\qquad\qquad\qquad\qquad
\end{equation*}
\begin{equation*}
\qquad\qquad\qquad\qquad\qquad\qquad\qquad\qquad\qquad\times
(\sum_{m=0}^{\infty}\frac{B_{m,\chi,\xi^{w_{1}w_{2}}}(w_{3}y_{3})}
{m!}(w_{1}w_{2}t)^{m})\\
\end{equation*}
\begin{equation}\label{a19}
=\sum_{n=0}^{\infty}(\sum_{k+l+m=n}^{}\binom{n}{k,l,m}
B_{k,\chi,\xi^{w_{2}w_{3}}}(w_{1}y_{1})
B_{l,\chi,\xi^{w_{1}w_{3}}}(w_{2}y_{2})
B_{m,\chi,\xi^{w_{1}w_{2}}}(w_{3}y_{3})\qquad\qquad
\end{equation}
\begin{equation*}
\qquad\qquad\qquad\qquad\qquad\qquad\qquad\qquad\qquad\times
w_{1}^{l+m}w_{2}^{k+m}w_{3}^{k+l})\frac{t^{n}}{n!},\qquad\qquad\\
\end{equation*}
where the inner sum is over all nonnegative integers $k,l,m$, with
$k+l+m=n$, and
\begin{equation*}
\binom{n}{k,l,m}=\frac{n!}{k!l!m!}.
\end{equation*}
\\
\\
(a-1) Here we write $I(\Lambda_{23}^{1})$ in two different ways:
\\
\\
(1)
$I(\Lambda_{23}^{1})$
\begin{equation}\label{a20}
\begin{split}
=\frac{1}{w_{3}}\int_{X}\chi(x_{1})
\xi^{w_{2}w_{3}x_{1}}e^{w_{2}w_{3}(x_{1}+w_{1}y_{1})t}&dx_{1}
\int_{X}\chi(x_{2})\xi^{w_{1}w_{3}x_{2}}
e^{w_{1}w_{3}(x_{2}+w_{2}y_{2})t}dx_{2}\\
&\times
\frac{dw_{3}\int_{X}\chi(x_{3})\xi^{w_{1}w_{2}x_{3}}
e^{w_{1}w_{2}x_{3}t}dx_{3}}{\int_{X}
\xi^{dw_{1}w_{2}w_{3}x_{4}}e^{dw_{1}w_{2}w_{3}x_{4}t}dx_{4}}\qquad
\end{split}
\end{equation}
\begin{equation*}
\begin{split}
=\frac{1}{w_{3}}(\sum_{k=0}^{\infty}B_{k,\chi,\xi^{w_{2}w_{3}}}
(w_{1}y_{1})\frac{(w_{2}w_{3}t)^{k}}{k!})
&(\sum_{l=0}^{\infty}B_{l,\chi,\xi^{w_{1}w_{3}}}(w_{2}y_{2})
\frac{(w_{1}w_{3}t)^{l}}{l!})\\
&\times
(\sum_{m=0}^{\infty}S_{m}(dw_{3}-1;\chi,\xi^{w_{1}w_{2}})
\frac{(w_{1}w_{2}t)^{m}}{m!})
\end{split}
\end{equation*}
\begin{equation}\label{a21}
\begin{split}
=\sum_{n=0}^{\infty}(\sum_{k+l+m=n}^{}\binom{n}{k,l,m}
&B_{k,\chi,\xi^{w_{2}w_{3}}}(w_{1}y_{1})
B_{l,\chi,\xi^{w_{1}w_{3}}}(w_{2}y_{2})\\
&\times S_{m}(dw_{3}-1;\chi,\xi^{w_{1}w_{2}})
w_{1}^{l+m}w_{2}^{k+m}w_{3}^{k+l-1})\frac{t^{n}}{n!}.\qquad
\end{split}
\end{equation}
\\
(2) Invoking (\ref{a9}), (\ref{a20}) can also be written as\\

$I(\Lambda_{23}^{1})$
\begin{equation*}
\begin{split}
=\frac{1}{w_{3}}\sum_{a=0}^{dw_{3}-1}\chi(a)\xi^{aw_{1}w_{2}}
&\int_{X}\chi(x_{1})\xi^{w_{2}w_{3}x_{1}}e^{w_{2}w_{3}
(x_{1}+w_{1}y_{1})t}dx_{1}\\
&\times\int_{X}\chi(x_{2})\xi^{w_{1}w_{3}x_{2}}
e^{w_{1}w_{3}(x_{2}+w_{2}y_{2}+\frac{w_{2}}{w_{3}}a)t}dx_{2}\qquad\\
\end{split}
\end{equation*}
\begin{equation*}
\begin{split}
=\frac{1}{w_{3}}\sum_{a=0}^{dw_{3}-1}\chi(a)\xi^{aw_{1}w_{2}}
&(\sum_{k=0}^{\infty}B_{k,\chi,\xi^{w_{2}w_{3}}}(w_{1}y_{1})
\frac{(w_{2}w_{3}t)^{k}}{k!})\\
&\times(\sum_{l=0}^{\infty}B_{l,\chi,\xi^{w_{1}w_{3}}}
(w_{2}y_{2}+\frac{w_{2}}{w_{3}}a)\frac{(w_{1}w_{3}t)^{l}}{l!})
\quad\qquad\\
\end{split}
\end{equation*}
\begin{equation}\label{a22}
\begin{split}
=\sum_{n=0}^{\infty}&(w_{3}^{n-1}\sum_{k=0}^{n}\binom{n}{k}
B_{k,\chi,\xi^{w_{2}w_{3}}}(w_{1}y_{1})\\
&\times\sum_{a=0}^{dw_{3}-1}\chi(a)\xi^{aw_{1}w_{2}}
B_{n-k,\chi,\xi^{w_{1}w_{3}}}(w_{2}y_{2}+\frac{w_{2}}{w_{3}}a)
w_{1}^{n-k}w_{2}^{k})\frac{t^{n}}{n!}.\qquad\qquad\qquad
\end{split}
\end{equation}
\\
(a-2) Here we write $I(\Lambda_{23}^{2})$ in three different ways:
\\
(1)
$I(\Lambda_{23}^{2})$
\begin{equation}\label{a23}
=\frac{1}{w_{2}w_{3}}\int_{X}\chi(x_{1})
\xi^{w_{2}w_{3}x_{1}}e^{w_{2}w_{3}(x_{1}+w_{1}y_{1})t}dx_{1}
\times
\frac{dw_{2}\int_{X}\chi(x_{2})\xi^{w_{1}w_{3}x_{2}}
e^{w_{1}w_{3}x_{2}t}dx_{2}}{\int_{X}
\xi^{dw_{1}w_{2}w_{3}x_{4}}
e^{dw_{1}w_{2}w_{3}x_{4}t}dx_{4}}\qquad\qquad\qquad\qquad\qquad
\end{equation}
\begin{equation*}
\qquad\qquad\qquad\qquad\times
\frac{dw_{3}\int_{X}\chi(x_{3})\xi^{w_{1}w_{2}x_{3}}e^{w_{1}w_{2}x_{3}t}
dx_{3}}{\int_{X}\xi^{dw_{1}w_{2}w_{3}x_{4}}
e^{dw_{1}w_{2}w_{3}x_{4}t}dx_{4}}\qquad
\end{equation*}
\begin{equation*}
=\frac{1}{w_{2}w_{3}}(\sum_{k=0}^{\infty}{B_{k,\chi,\xi^{w_{2}w_{3}}}
(w_{1}y_{1})\frac{(w_{2}w_{3}t)^k}{k!}})
(\sum_{l=0}^{\infty}{S_{l}(dw_{2}-1;\chi,\xi^{w_{1}w_{3}})
\frac{(w_{1}w_{3}t)^l}{l!}})\qquad
\end{equation*}
\begin{equation*}
\qquad\qquad\qquad\qquad\qquad\qquad\qquad\times
(\sum_{m=0}^{\infty}{S_{m}(dw_{3}-1;\chi,\xi^{w_{1}w_{2}})
\frac{(w_{1}w_{2}t)^m}{m!}})\quad\qquad
\end{equation*}
\begin{equation}\label{a24}
=\sum_{n=0}^{\infty}(\sum_{k+l+m=n}\binom{n}{k,l,m}B_{k,\chi,\xi^{w_{2}w_{3}}}
(w_{1}y_{1})S_{l}(dw_{2}-1;\chi,\xi^{w_{1}w_{3}})\qquad\qquad\qquad\qquad\qquad\qquad
\end{equation}
\begin{equation*}
\qquad\qquad\qquad\qquad\qquad\qquad\qquad\times
S_{m}(dw_{3}-1;\chi,\xi^{w_{1}w_{2}})w_{1}^{l+m}w_{2}^{k+m-1}
w_{3}^{k+l-1})\frac{t^{n}}{n!}.\qquad\qquad\qquad\qquad\qquad\qquad\qquad\qquad\qquad
\quad\quad
\end{equation*}
\\
\\
(2) Invoking (\ref{a9}), (\ref{a23}) can also be written as
\\
\\
~$I(\Lambda_{23}^{2})$
\begin{equation}\label{a25}
\begin{split}
=\frac{1}{w_{2}w_{3}}\sum_{a=0}^{dw_{2}-1}\chi(a)\xi^{aw_{1}w_{3}}
&\int_{X}\chi(x_{1})\xi^{w_{2}w_{3}x_{1}}e^{w_{2}w_{3}
(x_{1}+w_{1}y_{1}+\frac{w_{1}}{w_{2}}a)t}dx_{1}\\
&\times
\frac{dw_{3}\int_{X}\chi(x_{3})\xi^{w_{1}w_{2}x_{3}}e^{w_{1}w_{2}x_{3}t}
dx_{3}}{\int_{X}\xi^{dw_{1}w_{2}w_{3}x_{4}}e^{dw_{1}w_{2}w_{3}x_{4}t}dx_{4}}
\end{split}
\end{equation}
\begin{equation*}
\begin{split}
=\frac{1}{w_{2}w_{3}}\sum_{a=0}^{dw_{2}-1}\chi(a)
\xi^{aw_{1}w_{3}}&(\sum_{k=0}^{\infty}B_{k,\chi,\xi^{w_{2}w_{3}}}
(w_{1}y_{1}+\frac{w_{1}}{w_{2}}a)\frac{(w_{2}w_{3}t)^k}{k!})\\
&\times
(\sum_{l=0}^{\infty}S_{l}(dw_{3}-1;\chi,\xi^{w_{1}w_{2}})
\frac{(w_{1}w_{2}t)^l}{l!})
\end{split}
\end{equation*}
\begin{equation}\label{a26}
\begin{split}
=\sum_{n=0}^{\infty}(w_{2}^{n-1}\sum_{k=0}^{n}\binom{n}{k}
&\sum_{a=0}^{dw_{2}-1}\chi(a)\xi^{aw_{1}w_{3}}
B_{k,\chi,\xi^{w_{2}w_{3}}}(w_{1}y_{1}+\frac{w_{1}}{w_{2}}a)\\
&\times S_{n-k}(dw_{3}-1;\chi,\xi^{w_{1}w_{2}})w_{1}^{n-k}w_{3}^{k-1})\frac{t^{n}}{n!}.\qquad\qquad
\end{split}
\end{equation}
\\
\\
(3) Invoking (\ref{a9}) once again, (\ref{a25}) can be written as
\\
\\
\begin{equation*}
\begin{split}
I(\Lambda_{23}^{2})=\frac{1}{w_{2}w_{3}}\sum_{a=0}^{dw_{2}-1}\chi(a)
&\xi^{aw_{1}w_{3}}\sum_{b=0}^{dw_{3}-1}\chi(b)\xi^{aw_{1}w_{2}}\\
&\times\int_{X}\chi(x_{1})\xi^{w_{2}w_{3}x_{1}}e^{w_{2}w_{3}(x_{1}+w_{1}y_{1}+\frac{w_{1}}{w_{2}}a+\frac{w_{1}}{w_{3}}b)t}dx_{1}~
\end{split}
\end{equation*}
\begin{equation*}
\begin{split}
\qquad\qquad=\frac{1}{w_{2}w_{3}}\sum_{a=0}^{dw_{2}-1}
\chi(a)&\xi^{aw_{1}w_{3}}\sum_{b=0}^{dw_{3}-1}\chi(b)\xi^{bw_{1}w_{2}}\\
&\times(\sum_{n=0}^{\infty}B_{n,\chi,\xi^{w_{2}w_{3}}}
(w_{1}y_{1}+\frac{w_{1}}{w_{2}}a+\frac{w_{1}}{w_{3}}b)
\frac{(w_{2}w_{3}t)^n}{n!})
\end{split}
\end{equation*}
\begin{equation}\label{a27}
\begin{split}
=\sum_{n=0}^{\infty}((w_{2}w_{3})^{n-1}&\sum_{a=0}^{dw_{2}-1}
\sum_{b=0}^{dw_{3}-1}\chi(ab)\xi^{w_{1}(aw_{3}+bw_{2})}\\
&\times B_{n,\chi,\xi^{w_{2}w_{3}}}(w_{1}y_{1}+\frac{w_{1}}{w_{2}}a+\frac{w_{1}}{w_{3}}b))\frac{t^{n}}{n!}.\quad
\end{split}
\end{equation}
\\
\\
(a-3)
\\
\\
\begin{equation*}
I(\Lambda_{23}^{3})=\frac{1}{w_{1}w_{2}w_{3}}\times
\frac{dw_{1}\int_{X}\chi(x_{1})\xi^{w_{2}w_{3}x_{1}}e^{w_{2}w_{3}x_{1}t}
dx_{1}}{\int_{X}\xi^{dw_{1}w_{2}w_{3}x_{4}}e^{dw_{1}w_{2}w_{3}x_{4}t}
dx_{4}}\qquad\qquad\qquad\qquad\qquad\qquad\qquad
\end{equation*}
\begin{equation*}
\qquad\qquad\times\frac{dw_{2}\int_{X}\chi(x_{2})\xi^{w_{1}w_{3}x_{2}}
e^{w_{1}w_{3}x_{2}t}dx_{2}}{\int_{X}\xi^{dw_{1}w_{2}w_{3}x_{4}}
e^{dw_{1}w_{2}w_{3}x_{4}t}dx_{4}}
\times\frac{dw_{3}\int_{X}\chi(x_{3})
\xi^{w_{1}w_{2}x_{3}}e^{w_{1}w_{2}x_{3}t}dx_{3}}
{\int_{X}\xi^{dw_{1}w_{2}w_{3}x_{4}}e^{dw_{1}w_{2}w_{3}x_{4}t}dx_{4}}
\end{equation*}
\begin{align*}
=\frac{1}{w_{1}w_{2}w_{3}}(\sum_{k=0}^{\infty}S_{k}(dw_{1}-1;\chi,
\xi^{w_{2}w_{3}})&\frac{(w_{2}w_{3}t)^{k}}{k!})
(\sum_{l=0}^{\infty}S_{l}(dw_{2}-1;\chi,\xi^{w_{1}w_{3}})
\frac{(w_{1}w_{3}t)^{l}}{l!})\quad\qquad\\
&\times(\sum_{m=0}^{\infty}S_{m}(dw_{3}-1;\chi,\xi^{w_{1}w_{2}})
\frac{(w_{1}w_{2}t)^{m}}{m!})
\end{align*}
\begin{equation}\label{a28}
\begin{split}
=\sum_{n=0}^{\infty}(\sum_{k+l+m=n}^{}\binom{n}{k,l,m}&S_{k}
(dw_{1}-1;\chi,\xi^{w_{2}w_{3}})S_{l}(dw_{2}-1;\chi,\xi^{w_{1}w_{3}})\\
&\times S_{m}(dw_{3}-1;\chi,\xi^{w_{1}w_{2}})w_{1}^{l+m-1}w_{2}^{k+m-1}
w_{3}^{k+l-1})\frac{t^{n}}{n!}.\quad
\end{split}
\end{equation}
\\
(b) For Type $\Lambda_{13}^{i}~(i=0,1,2,3)$, we may consider the analogous things to the
ones in (a-0), (a-1), (a-2), and (a-3). However, each of those can be obtained from the corresponding ones in
(a-0), (a-1), (a-2), and (a-3). Indeed, if we substitute
    $w_{2}w_{3},w_{1}w_{3},w_{1}w_{2}$ respectively for $w_{1},w_{2},w_{3}$
in (\ref{a11}), this amounts to replacing $t$ by $w_{1}w_{2}w_{3}t$ and $\xi$ by $\xi^{w_{1}w_{2}w_{3}}$ in (\ref{a13}).
So, upon replacing  $w_{1},w_{2},w_{3}$ respectively by $w_{2}w_{3},w_{1}w_{3},w_{1}w_{2}$
, dividing by $(w_{1}w_{2}w_{3})^n$, and replacing
$\xi^{w_{1}w_{2}w_{3}}$ by $\xi$, in each of the expressions
of (\ref{a19}), (\ref{a21}), (\ref{a22}), (\ref{a24}), (\ref{a26})-(\ref{a28}), we will get the corresponding symmetric identities
for  Type $\Lambda_{13}^{i}~(i=0,1,2,3)$.
\\
\\
(c-0)
\begin{align*}
&I(\Lambda_{12}^{0})\\
&=\int_{X}\chi(x_{1})\xi^{w_{1}x_{1}}e^{w_{1}(x_{1}+w_{2}y)t}dx_{1}
\int_{X}\chi(x_{2})\xi^{w_{2}x_{2}}e^{w_{2}(x_{2}+w_{3}y)t}dx_{2}\\
&\qquad\qquad\qquad\qquad\qquad\qquad\qquad\qquad\qquad\times
\int_{X}\chi(x_{3})\xi^{w_{3}x_{3}}e^{w_{3}(x_{3}+w_{1}y)t}dx_{3}\qquad\qquad\qquad\\
&=(\sum_{n=0}^{\infty}\frac{B_{k,\chi,\xi^{w_{1}}}(w_{2}y)}{k!}
(w_{1}t)^{k})(\sum_{l=0}^{\infty}\frac{B_{l,\chi,\xi^{w_{2}}}
(w_{3}y)}{l!}(w_{2}t)^{l})
(\sum_{m=0}^{\infty}\frac{B_{m,\chi,\xi^{w_{3}}}(w_{1}y)}{m!}(w_{3}t)^{m})\qquad\qquad\qquad
\end{align*}
\begin{align}\label{a29}
=\sum_{n=0}^{\infty}(\sum_{k+l+m=n}\binom{n}{k,l,m}
B_{k,\chi,\xi^{w_{1}}}(w_{2}y)B_{l,\chi,\xi^{w_{2}}}(w_{3}y)
B_{m,\chi,\xi^{w_{3}}}(w_{1}y)w_{1}^{k}w_{2}^{l}w_{3}^{m})
\frac{t^n}{n!}.~\qquad\qquad
\end{align}
\\
\\
(c-1)
\begin{align*}
\begin{split}
I(\Lambda_{12}^{1})=\frac{1}{w_{1}w_{2}w_{3}}\frac{dw_{2}
\int_{X}\chi(x_{1})\xi^{w_{1}x_{1}}e^{w_{1}x_{1}t}dx_{1}}
{\int_{X}\xi^{dw_{1}w_{2}z_{3}}e^{dw_{1}w_{2}z_{3}t}dz_{3}}
&\times\frac{dw_{3}\int_{X}\chi(x_{2})\xi^{w_{2}x_{2}}
e^{w_{2}x_{2}t}dx_{2}}{\int_{X}\xi^{dw_{2}w_{3}z_{1}}
e^{dw_{2}w_{3}z_{1}t}dz_{1}}\\
&\times
\frac{dw_{1}\int_{X}\chi(x_{3})\xi^{w_{3}x_{3}}
e^{w_{3}x_{3}t}dx_{3}}{\int_{X}\xi^{dw_{3}w_{1}z_{2}}
e^{dw_{3}w_{1}z_{2}t}dz_{2}}
\end{split}
\end{align*}
\begin{align*}
\begin{split}
\qquad\qquad=\frac{1}{w_{1}w_{2}w_{3}}(\sum_{k=0}^{\infty}S_{k}
(dw_{2}-1;\chi,\xi^{w_{1}})&\frac{(w_{1}t)^{k}}{k!})
(\sum_{l=0}^{\infty}S_{l}(dw_{3}-1;\chi,\xi^{w_{2}})
\frac{(w_{2}t)^{l}}{l!})\\
&\times
(\sum_{m=0}^{\infty}S_{m}(dw_{1}-1;\chi,\xi^{w_{3}})
\frac{(w_{3}t)^{m}}{m!})
\end{split}
\end{align*}
\begin{equation}\label{a30}
\begin{split}
=\sum_{n=0}^{\infty}(\sum_{k+l+m=n}^{}\binom{n}{k,l,m}&S_{k}
(dw_{2}-1;\chi,\xi^{w_{1}})S_{l}(dw_{3}-1;\chi,\xi^{w_{2}})\qquad\qquad\\
&\times S_{m}(dw_{1}-1;\chi,\xi^{w_{3}})w_{1}^{k-1}w_{2}^{l-1}w_{3}^{m-1})
\frac{t^{n}}{n!}.
\end{split}
\end{equation}
\section{Main theorems}
As we noted earlier in the last paragraph of Section 2, the various
types of quotients of $p$-adic integrals are invariant under any
permutation of $w_{1},w_{2},w_{3}$. So the corresponding expressions
in Section 3 are also invariant under any permutation of
$w_{1},w_{2},w_{3}$. Thus our results about identities of symmetry
will be immediate consequences of this observation.

However, not all permutations of an expression in Section 3 yield
distinct ones. In fact, as these expressions are obtained by
permuting $w_{1},w_{2},w_{3}$ in a single one labeled by them, they
can be viewed as a group in a natural manner and hence it is
isomorphic to a quotient of $S_{3}$. In particular, the number of
possible distinct expressions are $1,2,3$ or $6$. (a-0), (a-1(1)),
(a-1(2)), and (a-2(2)) give the full six identities of symmetry,
(a-2(1)) and (a-2(3)) yield three identities of symmetry, and (c-0)
and (c-1) give two identities of symmetry, while the expression in
(a-3) yields no identities of symmetry.

Here we will just consider the cases of Theorems \ref{A4} and \ref{A8}, leaving the others as easy
exercises for the reader. As for the case of Theorem \ref{A4}, in addition to (\ref{a41})-(\ref{a43}),
we get the following three ones:\\

\begin{align}
\begin{split}
\label{a31}
\sum_{k+l+m=n}\binom{n}{k,l,m}&B_{k,\chi,\xi^{w_{2}w_{3}}}(w_{1}y_{1})
S_{l}(dw_{3}-1;\chi,\xi^{w_{1}w_{2}})\\
&\times S_{m}(dw_{2}-1;\chi,\xi^{w_{1}w_{3}})w_{1}^{l+m}w_{3}^{k+m-1}
w_{2}^{k+l-1},
\end{split}
\end{align}
\begin{align}
\begin{split}
\label{a32}
\sum_{k+l+m=n}\binom{n}{k,l,m}&B_{k,\chi,\xi^{w_{1}w_{3}}}(w_{2}y_{1})
S_{l}(dw_{1}-1;\chi,\xi^{w_{2}w_{3}})\\
&\times S_{m}(dw_{3}-1;\chi,\xi^{w_{1}w_{2}})
w_{2}^{l+m}w_{1}^{k+m-1}w_{3}^{k+l-1},
\end{split}
\end{align}
\begin{align}
\begin{split}
\label{a33}
\sum_{k+l+m=n}\binom{n}{k,l,m}&B_{k,\chi,\xi^{w_{1}w_{2}}}(w_{3}y_{1})
S_{l}(dw_{2}-1;\chi,\xi^{w_{1}w_{3}})\\
&\times S_{m}(dw_{1}-1;\chi,\xi^{w_{2}w_{3}})
w_{3}^{l+m}w_{2}^{k+m-1}w_{1}^{k+l-1}.
\end{split}
\end{align}
But, by interchanging $l$ and $m$, we see that (\ref{a31}),
(\ref{a32}), and (\ref{a33}) are respectively equal to (\ref{a41}),
(\ref{a42}), and (\ref{a43}).
\\
As to Theorem \ref{A8}, in addition to (\ref{a47}) and (\ref{a48}), we have:

\begin{align}
\begin{split}
\label{a34}
\sum_{k+l+m=n}\binom{n}{k,l,m}&S_{k}(dw_{2}-1;\chi,\xi^{w_{1}})
S_{l}(dw_{3}-1;\chi,\xi^{w_{2}})\\
&\times S_{m}(dw_{1}-1;\chi,\xi^{w_{3}})w_{1}^{k-1}w_{2}^{l-1}w_{3}^{m-1},
\end{split}
\end{align}
\begin{align}
\begin{split}
\label{a35}
\sum_{k+l+m=n}\binom{n}{k,l,m}&S_{k}(dw_{3}-1;\chi,\xi^{w_{2}})
S_{l}(dw_{1}-1;\chi,\xi^{w_{3}})\\
&\times S_{m}(dw_{2}-1;\chi,\xi^{w_{1}})w_{2}^{k-1}w_{3}^{l-1}w_{1}^{m-1},
\end{split}
\end{align}
\begin{align}
\begin{split}
\label{a36}
\sum_{k+l+m=n}\binom{n}{k,l,m}&S_{k}(dw_{3}-1;\chi,\xi^{w_{1}})
S_{l}(dw_{2}-1;\chi,\xi^{w_{3}})\\
&\times S_{m}(dw_{1}-1;\chi,\xi^{w_{2}})w_{1}^{k-1}w_{3}^{l-1}w_{2}^{m-1},
\end{split}
\end{align}
\begin{align}
\begin{split}
\label{a37}
\sum_{k+l+m=n}\binom{n}{k,l,m}&S_{k}(dw_{2}-1;\chi,\xi^{w_{3}})
S_{l}(dw_{1}-1;\chi,\xi^{w_{2}})\\
&\times S_{m}(dw_{3}-1;\chi,\xi^{w_{1}})w_{3}^{k-1}w_{2}^{l-1}w_{1}^{m-1}.
\end{split}
\end{align}
\\
  However, (\ref{a34}) and (\ref{a35}) are equal to (\ref{a47}), as we can see by applying the permutations $k\rightarrow l,l\rightarrow
m,m\rightarrow k$ for (\ref{a34}) and $k\rightarrow m,l\rightarrow
k,m\rightarrow l$ for (\ref{a35}). Similarly, we see that (\ref{a36})
and (\ref{a37}) are equal to (\ref{a48}), by applying permutations
$k\rightarrow l,l\rightarrow m,m\rightarrow k$ for (\ref{a36}) and
$k\rightarrow m,l\rightarrow k,m\rightarrow l$ for (\ref{a37}).
\begin{theorem}\label{A1}
Let $w_{1},w_{2},w_{3}$ be any positive integers. Then we have:
\begin{equation}\label{a38}
\begin{split}
\sum_{k+l+m=n}\binom{n}{k,l,m}B_{k,\chi,\xi^{w_{2}w_{3}}}
(w_{1}y_{1})&B_{l,\chi,\xi^{w_{1}w_{3}}}(w_{2}y_{2})\\
&\times B_{m,\chi,\xi^{w_{1}w_{2}}}(w_{3}y_{3})w_{1}^{l+m}w_{2}^{k+m}w_{3}^{k+l}
\end{split}
\end{equation}
\begin{equation*}
\begin{split}
=\sum_{k+l+m=n}\binom{n}{k,l,m}B_{k,\chi,\xi^{w_{2}w_{3}}}
(w_{1}y_{1})&B_{l,\chi,\xi^{w_{1}w_{2}}}(w_{3}y_{2})\\
&\times B_{m,\chi,\xi^{w_{1}w_{3}}}(w_{2}y_{3})w_{1}^{l+m}w_{3}^{k+m}w_{2}^{k+l}
\end{split}
\end{equation*}
\begin{equation*}
\begin{split}
=\sum_{k+l+m=n}\binom{n}{k,l,m}B_{k,\chi,\xi^{w_{1}w_{3}}}
(w_{2}y_{1})&B_{l,\chi,\xi^{w_{2}w_{3}}}(w_{1}y_{2})\\
&\times B_{m,\chi,\xi^{w_{1}w_{2}}}(w_{3}y_{3})w_{2}^{l+m}w_{1}^{k+m}w_{3}^{k+l}
\end{split}
\end{equation*}
\begin{equation*}
\begin{split}
=\sum_{k+l+m=n}\binom{n}{k,l,m}B_{k,\chi,\xi^{w_{1}w_{3}}}
(w_{2}y_{1})&B_{l,\chi,\xi^{w_{1}w_{2}}}(w_{3}y_{2})\\
&\times B_{m,\chi,\xi^{w_{2}w_{3}}}(w_{1}y_{3})w_{2}^{l+m}w_{3}^{k+m}w_{1}^{k+l}
\end{split}
\end{equation*}
\begin{equation*}
\begin{split}
=\sum_{k+l+m=n}\binom{n}{k,l,m}B_{k,\chi,\xi^{w_{1}w_{2}}}
(w_{3}y_{1})&B_{l,\chi,\xi^{w_{2}w_{3}}}(w_{1}y_{2})\\
&\times B_{m,\chi,\xi^{w_{1}w_{3}}}(w_{2}y_{3})w_{3}^{l+m}w_{1}^{k+m}w_{2}^{k+l}
\end{split}
\end{equation*}
\begin{equation*}
\begin{split}
=\sum_{k+l+m=n}\binom{n}{k,l,m}B_{k,\chi,\xi^{w_{1}w_{2}}}
(w_{3}y_{1})&B_{l,\chi,\xi^{w_{1}w_{3}}}(w_{2}y_{2})\\
&\times B_{m,\chi,\xi^{w_{2}w_{3}}}(w_{1}y_{3})w_{3}^{l+m}w_{2}^{k+m}w_{1}^{k+l}.
\end{split}
\end{equation*}
\end{theorem}
\begin{theorem}\label{A2}
  Let $w_{1},w_{2},w_{3}$ be any positive integers. Then we have:
\begin{equation*}
\begin{split}
\sum_{k+l+m=n}\binom{n}{k,l,m}&B_{k,\chi,\xi^{w_{2}w_{3}}}
(w_{1}y_{1})B_{l,\chi,\xi^{w_{1}w_{3}}}(w_{2}y_{2})\\
&\times S_{m}(dw_{3}-1;\chi,\xi^{w_{1}w_{2}})w_{1}^{l+m}w_{2}^{k+m}w_{3}^{k+l-1}\\
\end{split}
\end{equation*}
\begin{equation*}
\begin{split}
=\sum_{k+l+m=n}\binom{n}{k,l,m}&B_{k,\chi,\xi^{w_{2}w_{3}}}
(w_{1}y_{1})B_{l,\chi,\xi^{w_{1}w_{2}}}(w_{3}y_{2})\\
&\times S_{m}(dw_{2}-1;\chi,\xi^{w_{1}w_{3}})w_{1}^{l+m}w_{3}^{k+m}w_{2}^{k+l-1}\\
\end{split}
\end{equation*}
\begin{equation}\label{a39}
\begin{split}
=\sum_{k+l+m=n}\binom{n}{k,l,m}&B_{k,\chi,\xi^{w_{1}w_{3}}}
(w_{2}y_{1})B_{l,\chi,\xi^{w_{2}w_{3}}}(w_{1}y_{2})\\
&\times S_{m}(dw_{3}-1;\chi,\xi^{w_{1}w_{2}})w_{2}^{l+m}w_{1}^{k+m}w_{3}^{k+l-1}\qquad\\
\end{split}
\end{equation}
\begin{equation*}
\begin{split}
=\sum_{k+l+m=n}\binom{n}{k,l,m}&B_{k,\chi,\xi^{w_{1}w_{3}}}
(w_{2}y_{1})B_{l,\chi,\xi^{w_{1}w_{2}}}(w_{3}y_{2})\\
&\times S_{m}(dw_{1}-1;\chi,\xi^{w_{2}w_{3}})w_{2}^{l+m}w_{3}^{k+m}w_{1}^{k+l-1}\\
\end{split}
\end{equation*}
\begin{equation*}
\begin{split}
=\sum_{k+l+m=n}\binom{n}{k,l,m}&B_{k,\chi,\xi^{w_{1}w_{2}}}(w_{3}y_{1})
B_{l,\chi,\xi^{w_{1}w_{3}}}(w_{2}y_{2})\\
&\times S_{m}(dw_{1}-1;\chi,\xi^{w_{2}w_{3}})w_{3}^{l+m}w_{2}^{k+m}w_{1}^{k+l-1}\\
\end{split}
\end{equation*}
\begin{equation*}
\begin{split}
=\sum_{k+l+m=n}\binom{n}{k,l,m}&B_{k,\chi,\xi^{w_{1}w_{2}}}
(w_{3}y_{1})B_{l,\chi,\xi^{w_{2}w_{3}}}(w_{1}y_{2})\\
&\times S_{m}(dw_{2}-1;\chi,\xi^{w_{1}w_{3}})w_{3}^{l+m}w_{1}^{k+m}w_{2}^{k+l-1}.
\end{split}
\end{equation*}
\end{theorem}

\begin{theorem}\label{A3}
Let $w_{1},w_{2},w_{3}$ be any positive integers. Then we have:
\begin{equation}\label{a40}
\begin{split}
&w_{1}^{n-1}\sum_{k=0}^{n}\binom{n}{k}B_{k,\chi,\xi^{w_{1}w_{2}}}
(w_{3}y_{1})\sum_{a=0}^{dw_{1}-1}\chi(a)\xi^{aw_{2}w_{3}}
B_{n-k,\chi,\xi^{w_{1}w_{3}}}(w_{2}y_{2}+\frac{w_{2}}{w_{1}}a)w_{3}^{n-k}w_{2}^{k}\\
=&w_{1}^{n-1}\sum_{k=0}^{n}\binom{n}{k}B_{k,\chi,\xi^{w_{1}w_{3}}}
(w_{2}y_{1})\sum_{a=0}^{dw_{1}-1}\chi(a)\xi^{aw_{2}w_{3}}
B_{n-k,\chi,\xi^{w_{1}w_{2}}}(w_{3}y_{2}+\frac{w_{3}}{w_{1}}a)w_{2}^{n-k}w_{3}^{k}\\
=&w_{2}^{n-1}\sum_{k=0}^{n}\binom{n}{k}B_{k,\chi,\xi^{w_{1}w_{2}}}
(w_{3}y_{1})\sum_{a=0}^{dw_{2}-1}\chi(a)\xi^{aw_{1}w_{3}}
B_{n-k,\chi,\xi^{w_{2}w_{3}}}(w_{1}y_{2}+\frac{w_{1}}{w_{2}}a)w_{3}^{n-k}w_{1}^{k}\\
=&w_{2}^{n-1}\sum_{k=0}^{n}\binom{n}{k}B_{k,\chi,\xi^{w_{2}w_{3}}}
(w_{1}y_{1})\sum_{a=0}^{dw_{2}-1}\chi(a)\xi^{aw_{1}w_{3}}
B_{n-k,\chi,\xi^{w_{1}w_{2}}}(w_{3}y_{2}+\frac{w_{3}}{w_{2}}a)w_{1}^{n-k}w_{3}^{k}\\
=&w_{3}^{n-1}\sum_{k=0}^{n}\binom{n}{k}B_{k,\chi,\xi^{w_{1}w_{3}}}
(w_{2}y_{1})\sum_{a=0}^{dw_{3}-1}\chi(a)\xi^{aw_{1}w_{2}}
B_{n-k,\chi,\xi^{w_{2}w_{3}}}(w_{1}y_{2}+\frac{w_{1}}{w_{2}}a)w_{2}^{n-k}w_{1}^{k}\\
=&w_{3}^{n-1}\sum_{k=0}^{n}\binom{n}{k}B_{k,\chi,\xi^{w_{2}w_{3}}}
(w_{1}y_{1})\sum_{a=0}^{dw_{3}-1}\chi(a)\xi^{aw_{1}w_{2}}
B_{n-k,\chi,\xi^{w_{1}w_{3}}}(w_{2}y_{2}+\frac{w_{2}}{w_{3}}a)w_{1}^{n-k}w_{2}^{k}.
\end{split}
\end{equation}
\end{theorem}

\begin{theorem}\label{A4}
Let $w_{1},w_{2},w_{3}$ be any positive integers. Then we have the
following three symmetries in $w_{1},w_{2},w_{3}$:

\begin{equation}\label{a41}
\begin{split}
\sum_{k+l+m=n}\binom{n}{k,l,m}&B_{k,\chi,\xi^{w_{2}w_{3}}}
(w_{1}y_{1})S_{l}(dw_{2}-1;\chi,\xi^{w_{1}w_{3}})\\
&\times S_{m}(dw_{3}-1;\chi,\xi^{w_{1}w_{2}})w_{1}^{l+m}w_{2}^{k+m-1}w_{3}^{k+l-1}
\end{split}
\end{equation}
\begin{equation}
\begin{split}
\label{a42}
=\sum_{k+l+m=n}\binom{n}{k,l,m}&B_{k,\chi,\xi^{w_{1}w_{3}}}
(w_{2}y_{1})S_{l}(dw_{3}-1;\chi,\xi^{w_{1}w_{2}})\\
&\times S_{m}(dw_{1}-1;\chi,\xi^{w_{2}w_{3}})w_{2}^{l+m}w_{3}^{k+m-1}w_{1}^{k+l-1}
\end{split}
\end{equation}
\begin{equation}
\begin{split}
\label{a43}
=\sum_{k+l+m=n}\binom{n}{k,l,m}&B_{k,\chi,\xi^{w_{1}w_{2}}}
(w_{3}y_{1})S_{l}(dw_{1}-1,\chi,\xi^{w_{2}w_{3}})\\
&\times S_{m}(dw_{2}-1;\chi,\xi^{w_{1}w_{3}})w_{3}^{l+m}w_{1}^{k+m-1}w_{2}^{k+l-1}.
\end{split}
\end{equation}
\end{theorem}

\begin{theorem}\label{A5}
Let $w_{1},w_{2},w_{3}$ be any positive integers. Then we have:
\begin{equation}\label{a44}
\begin{split}
w_{1}^{n-1}\sum_{k=0}^{n}\binom{n}{k}\sum_{a=0}^{dw_{1}-1}
\chi(a)\xi^{aw_{2}w_{3}}&B_{k,\chi,\xi^{w_{1}w_{3}}}
(w_{2}y_{1}+\frac{w_{2}}{w_{1}}a)\\
&\times S_{n-k}(dw_{3}-1;\chi,\xi^{w_{1}w_{2}})
w_{2}^{n-k}w_{3}^{k-1}\\
\end{split}
\end{equation}
\begin{equation*}
\begin{split}
=w_{1}^{n-1}\sum_{k=0}^{n}\binom{n}{k}\sum_{a=0}^{dw_{1}-1}\chi(a)
\xi^{aw_{2}w_{3}}
&B_{k,\chi,\xi^{w_{1}w_{2}}}(w_{3}y_{1}+\frac{w_{3}}{w_{1}}a)\\
&\times S_{n-k}(dw_{2}-1;\chi,\xi^{w_{1}w_{3}})w_{3}^{n-k}w_{2}^{k-1}
\end{split}
\end{equation*}
\begin{equation*}
\begin{split}
=w_{2}^{n-1}\sum_{k=0}^{n}\binom{n}{k}\sum_{a=0}^{dw_{2}-1}\chi(a)
\xi^{aw_{1}w_{3}}&B_{k,\chi,\xi^{w_{2}w_{3}}}(w_{1}y_{1}+\frac{w_{1}}
{w_{2}}a)\\
&\times S_{n-k}(dw_{3}-1;\chi,\xi^{w_{1}w_{2}})w_{1}^{n-k}w_{3}^{k-1}
\end{split}
\end{equation*}
\begin{equation*}
\begin{split}
=w_{2}^{n-1}\sum_{k=0}^{n}\binom{n}{k}\sum_{a=0}^{dw_{2}-1}\chi(a)
\xi^{aw_{1}w_{3}}&B_{k,\chi,\xi^{w_{1}w_{2}}}(w_{3}y_{1}+\frac{w_{3}}
{w_{2}}a)\\
&\times S_{n-k}(dw_{1}-1;\chi,\xi^{w_{2}w_{3}})w_{3}^{n-k}w_{1}^{k-1}
\end{split}
\end{equation*}
\begin{equation*}
\begin{split}
=w_{3}^{n-1}\sum_{k=0}^{n}\binom{n}{k}\sum_{a=0}^{dw_{3}-1}\chi(a)
\xi^{aw_{1}w_{2}}&B_{k,\chi,\xi^{w_{2}w_{3}}}(w_{1}y_{1}+\frac{w_{1}}
{w_{3}}a)\\
&\times S_{n-k}(dw_{2}-1;\chi,\xi^{w_{1}w_{3}})w_{1}^{n-k}w_{2}^{k-1}
\end{split}
\end{equation*}
\begin{equation*}
\begin{split}
=w_{3}^{n-1}\sum_{k=0}^{n}\binom{n}{k}\sum_{a=0}^{dw_{3}-1}\chi(a)
\xi^{aw_{1}w_{2}}&B_{k,\chi,\xi^{w_{1}w_{3}}}(w_{2}y_{1}+\frac{w_{2}}
{w_{3}}a)\\
&\times S_{n-k}(dw_{1}-1;\chi,\xi^{w_{2}w_{3}})w_{2}^{n-k}w_{1}^{k-1}.
\end{split}
\end{equation*}
\end{theorem}

\begin{theorem}\label{A6}
Let $w_{1},w_{2},w_{3}$ be any positive integers. Then we have the
following three symmetries in $w_{1},w_{2},w_{3}$:\\
\begin{equation}\label{s1}
\begin{split}
&(w_{1}w_{2})^{n-1}\sum_{a=0}^{dw_{1}-1}\sum_{b=0}^{dw_{2}-1}
\chi(ab)\xi^{w_{3}(aw_{2}+bw_{1})}B_{n,\chi,\xi^{w_{1}w_{2}}}
(w_{3}y_{1}+\frac{w_{3}}{w_{1}}a+\frac{w_{3}}{w_{2}}b)\\
=&(w_{2}w_{3})^{n-1}\sum_{a=0}^{dw_{2}-1}\sum_{b=0}^{dw_{3}-1}
\chi(ab)\xi^{w_{1}(aw_{3}+bw_{2})}B_{n,\chi,\xi^{w_{2}w_{3}}}
(w_{1}y_{1}+\frac{w_{1}}{w_{2}}a+\frac{w_{1}}{w_{3}}b)\\
=&(w_{3}w_{1})^{n-1}\sum_{a=0}^{dw_{3}-1}\sum_{b=0}^{dw_{1}-1}
\chi(ab)\xi^{w_{2}(aw_{1}+bw_{3})}B_{n,\chi,\xi^{w_{1}w_{3}}}
(w_{2}y_{1}+\frac{w_{2}}{w_{3}}a+\frac{w_{2}}{w_{1}}b).
\end{split}
\end{equation}
\end{theorem}

\begin{theorem}\label{A7}
Let $w_{1},w_{2},w_{3}$ be any positive integers. Then we have the
following two symmetries in $w_{1},w_{2},w_{3}$:\\
\begin{equation}\label{a46}
\begin{split}
&\sum_{k+l+m=n}\binom{n}{k,l,m}B_{k,\chi,\xi^{w_{3}}}(w_{1}y)
B_{l,\chi,\xi^{w_{1}}}(w_{2}y)B_{m,\chi,\xi^{w_{2}}}(w_{3}y)w_{3}^{k}w_{1}^{l}w_{2}^{m}\\
=&\sum_{k+l+m=n}\binom{n}{k,l,m}B_{k,\chi,\xi^{w_{2}}}(w_{1}y)
B_{l,\chi,\xi^{w_{1}}}(w_{3}y)B_{m,\chi,\xi^{w_{3}}}(w_{2}y)w_{2}^{k}w_{1}^{l}w_{3}^{m}.
\end{split}
\end{equation}
\end{theorem}

\begin{theorem}\label{A8}
Let $w_{1},w_{2},w_{3}$ be any positive integers. Then we have
the following two symmetries in $w_{1},w_{2},w_{3}$:
\begin{align}
\begin{split}
\label{a47}
\sum_{k+l+m=n}\binom{n}{k,l,m}&S_{k}(dw_{1}-1;\chi,\xi^{w_{3}})
S_{l}(dw_{2}-1;\chi,\xi^{w_{1}})\\
&\times S_{m}(dw_{3}-1;\chi,\xi^{w_{2}})w_{3}^{k-1}w_{1}^{l-1}w_{2}^{m-1}\\
\end{split}
\end{align}
\begin{align}
\begin{split}
\label{a48}
=\sum_{k+l+m=n}\binom{n}{k,l,m}&S_{k}(dw_{1}-1;\chi,\xi^{w_{2}})
S_{l}(dw_{3}-1;\chi,\xi^{w_{1}})\\
&\times S_{m}(dw_{2}-1;\chi,\xi^{w_{3}})w_{2}^{k-1}w_{1}^{l-1}w_{3}^{m-1}.
\end{split}
\end{align}
\end{theorem}

\bibliographystyle{amsplain}

\end{document}